\documentclass[12pt]{article}
\usepackage{amsmath}
\usepackage{amssymb}
\usepackage{amsmath,amssymb,amsbsy,amsfonts,amsthm,latexsym,
amsopn,amstext,amsxtra,euscript,amscd}

\newtheorem{thm}{Theorem}
\newtheorem{theorem}[thm]{Theorem}

\begin{document}

\baselineskip 15pt

\title{\sc the sum-product estimate for large subsets of prime fields}

\author{{\sc M. Z. Garaev}
\\
\normalsize{Instituto de Matem{\'a}ticas}
\\
\normalsize{Universidad Nacional Aut\'onoma de M\'exico}
\\
\normalsize{Campus Morelia, Apartado Postal 61-3 (Xangari)}
\\
\normalsize{C.P. 58089, Morelia, Michoac{\'a}n, M{\'e}xico} \\
\normalsize{\tt garaev@matmor.unam.mx}\\
}

\date{\empty}

\pagenumbering{arabic}

\maketitle

\begin{abstract}
Let $\mathbb{F}_p$ be the field of a prime order $p.$ It is known
that for any integer $N\in [1,p]$ one can construct a subset
$A\subset\mathbb{F}_p$ with $|A|= N$ such that
$$
\max\{|A+A|, |AA|\}\ll p^{1/2}|A|^{1/2}.
$$
In the present paper we prove that if $A\subset \mathbb{F}_p$ with
$|A|>p^{2/3},$ then
$$
\max\{|A+A|, |AA|\}\gg p^{1/2}|A|^{1/2}.
$$
\end{abstract}



\section{Introduction}

Let $\mathbb{F}_p$ be the field of residue classes modulo a prime
number $p$ and let $A\subset \mathbb{F}_p.$ Consider the sum set
$$
A+A=\{a+b: \, a\in A, \, b\in A\}
$$
and the product set
$$
AA =\{ab: \, a\in A, \, b\in A\}.
$$
From the work of Bourgain, Katz, Tao~\cite{BKT} and Bourgain,
Glibichuk, Konyagin~\cite{BGK} it is known that if
$|A|<p^{1-\delta},$ where $\delta>0,$ then one has the sum-product
estimate
\begin{equation}
\label{eqn:sumprodclass} \max\{|A+A|,|AA|\}\gg |A|^{1+\varepsilon}
\end{equation}
for some $\varepsilon=\varepsilon(\delta)>0.$ This result and its
versions have found many important applications in various areas of
mathematics.

In the corresponding problem for integers (i.e., if the field
$\mathbb{F}_p$ is replaced by the set of integers) the conjecture of
Erd\"{o}s and Szemer\'edi~\cite{ErSz} is that $\max\{|A+A|,|AA|\}\gg
|A|^{2-\varepsilon}$ for any given $\varepsilon>0.$ At present the
best known bound in the integer problem is $\max\{|A+A|,|AA|\}\gg
|A|^{14/11}(\log|A|)^{-3/11}$ due to Solymosi~\cite{Sol}.

Explicit versions of~\eqref{eqn:sumprodclass} have been obtained
in~\cite{Gar}--\cite{KS1}. For subsets with relatively small
cardinalities (say, $|A|<p^{13/25}$), in~\cite{Gar} we proved that
$$
\max\{|A+A|,|AA|\}\gg |A|^{15/14}(\log|A|)^{O(1)}
$$
which was subsequently improved in~\cite{KS} to
$$
\max\{|A+A|,|AA|\}\gg |A|^{14/13}(\log|A|)^{O(1)}.
$$
In~\cite{Gar} we have also considered the case of subsets with
larger cardinalities, which had been previously studied
in~\cite{HIS}. We have shown, for example, that
$$
\max\{|A+A|, |AA|\}\gg \min\Bigl\{|A|^{2/3}p^{1/3},\,
|A|^{5/3}p^{-1/3}\Bigr\}(\log |A|)^{O(1)}.
$$
One may conjecture that the estimate
$$
\max\{|A+A|,|AA|\}\gg\min\{|A|^{2-\varepsilon},
|A|^{1/2}p^{1/2-\varepsilon}\}
$$
holds for all subsets $A\subset\mathbb{F}_p.$ The motivation for the
quantity $|A|^{1/2}p^{1/2-\varepsilon}$ is clear from the
construction in~\cite{MC} which can be described as follows. Let $g$
be a generator of $\mathbb{F}_p^*.$ By the pigeon-hole principle,
for any $N\in [1,p]$ and for any integer $M\approx p^{1/2}N^{1/2}$
(which we associate with $M\pmod p$),  there exists $L$ such that
$$
|\{g^{x}: 1\le x\le M\}\cap \{L+1, L+2, \ldots, L+M\}|\gg M^2/p\gg
N.
$$
Obviously, any subset $ A\subset \{g^{x}: 1\le x\le M\}\cap \{L+1,
L+2, \ldots, L+M\} $ with $|A|\approx N$ satisfies $\max\{|A+A|,
|AA|\}\ll p^{1/2}|A|^{1/2}.$

Thus, it follows that for any integer $N\in [1,p]$ there exists a
subset $A\subset \mathbb{F}_p$ with $|A|=N$ such that
$$
\max\{|A+A|, |AA|\}\ll p^{1/2}|A|^{1/2}.
$$
In the present paper we prove the following statement.

\begin{theorem}
\label{thm:sumprodestPrime} Let $A\subset \mathbb{F}_p.$ Then
$$
|A+A||AA|\gg \min\Bigl\{ p|A|,\, \frac{|A|^{4}}{p}\Bigr\}.
$$
\end{theorem}

In view of the foregoing discussion, in the range $|A|>p^{2/3}$ our
result implies the optimal in the general setting bound
$$
\max\{|A+A|, |AA|\}\gg p^{1/2}|A|^{1/2}.
$$

The following generalization of Theorem~\ref{thm:sumprodestPrime}
improves the corresponding result from~\cite{Vu}, where an analogy
of the sum-product estimate from~\cite{HIS} has been obtained for
subsets of $\mathbb{Z}_m,$ the ring of residue classes modulo $m.$
\begin{theorem}
\label{thm:sumprodestCom} Let $A\subset \mathbb{Z}_m.$  Then
$$
|A+A||AA|\gg \min\Bigl\{ m|A|,\, \frac{|A|^{4}}{m}\Bigl(\sum_{\substack {d|m\\
d<m}}d^{1/2}\Bigr)^{-2}\Bigr\}.
$$
\end{theorem}

We remark that if $m=p^2,\, A=\{p\, x: \, x\in \mathbb{Z}_m\},$
where $p$ is a prime number, then $|A|=|A+A|=m^{1/2},\, |AA|=1$ and
the left hand side of the estimate of
Theorem~\ref{thm:sumprodestCom} is of the same order of magnitude as
the right hand side.

\section{Proof of Theorem~\ref{thm:sumprodestPrime}}

We can assume that $\{0\}\not\in A.$ Consider the equation
\begin{equation}
\label{eqn:aa/a+a=a+a} xa_1^{-1}+a_2=y, \quad (x,\, a_1,\, a_2,\,
y)\in (AA)\times A\times A\times(A+A).
\end{equation}
For any triple $(a_1,\, a_2,\, a_3)\in A\times A\times A$ the vector
$$
(a_1a_3,\, a_1, \, a_2,\, a_3+a_2)\in (AA)\times A\times
A\times(A+A)
$$
is a solution of~\eqref{eqn:aa/a+a=a+a}. For different triples
$(a_1,\, a_2,\, a_3)\in A\times A\times A$ correspond different
solutions $(a_1a_3,\, a_1,\, a_2,\, a_3+a_2).$ Thus, the number $J$
of solutions of the equation~\eqref{eqn:aa/a+a=a+a} satisfies $J\ge
|A|^3.$ Expressing $J$ via additive characters and following the
standard procedure, we obtain
\begin{eqnarray*}
 |A|^3\le J=\frac{1}{p}\sum_{n=0}^{p-1}\sum_{x\in AA}\sum_{a_1\in
A}\sum_{a_2\in A}\sum_{y\in A+A}e_p(n(xa_1^{-1}+ a_2-y))\le \\
\frac{|AA||A|^2|A+A|}{p}
+\frac{1}{p}\sum_{n=1}^{p-1}\Bigl|\sum_{x\in AA}\sum_{a_1\in
A}e_p(nxa_1^{-1})\Bigr|\Bigl|\sum_{a_2\in A}\sum_{y\in
A+A}e_p(n(a_2-y))\Bigr|.
\end{eqnarray*}
Since
$$
\max_{(n,p)=1}\Bigl|\sum_{x\in AA}\sum_{a_1\in
A}e_p(nxa_1^{-1})\Bigr|\le \sqrt{p|AA||A|},
$$
(see, for example,~\cite[Chapter VI]{Vi}), we have
$$
|A|^3\le \frac{|AA||A|^2|A+A|}{p} +
\frac{\sqrt{p|AA||A|}}{p}\sum_{n=0}^{p-1}\Bigl|\sum_{a_2\in
A}e_p(na_2)\Bigr|\Bigl|\sum_{y\in A+A}e_p(ny)\Bigr|.
$$
Applying the Cauchy-Schwarz inequality to the sum over $n,$ we get
$$
|A|^3\le \frac{|AA||A|^2|A+A|}{p}+\sqrt{p|AA||A|}\sqrt{|A||A+A|}
$$
and the result follows.

\section{Proof of Theorem~\ref{thm:sumprodestCom}}

Let us first reduce the problem to the case $A\subset
\mathbb{Z}_m^*.$ Let
$$
d_0=\min\{\,(a,m)\,:\quad a\in A\}.
$$
Then,
$$
|AA|\ge |d_0A|\ge \frac{|A|}{d_0}.
$$
If
$$
|A|^2 \le \frac{4m}{d_0}\Bigl(\sum_{\substack {d|m\\
d<m}}d^{1/2}\Bigr)^{2},
$$
then
$$
\frac{|A|^{4}}{m}\Bigl(\sum_{\substack {d|m\\
d<m}}d^{1/2}\Bigr)^{-2} \le \frac{4|A|^2}{d_0}\le 4|A||AA|\le
4|A+A||AA|
$$
and the statement becomes trivial in this case. Thus, we can assume
that
$$
|A|^2 > \frac{4m}{d_0}\Bigl(\sum_{\substack {d|m\\
d<m}}d^{1/2}\Bigr)^{2}.
$$
But then we have that
\begin{eqnarray*}
&& \{a\in A:\,\, (a,m)>1\}=\{a\in A:\,\, (a,m)\ge \max\{d_0, 2\}\} \le  \sum_{\substack{d|m\\
d\ge \max\{d_0,2\}}}\frac{m}{d}\\ && =
\sum_{\substack{d|m\\
d\le \min\{m/d_0,m/2\}}}d  \,\, \le \,\,  \Bigl(\frac{m}{d_0}\Bigr)^{1/2}\sum_{\substack {d|m\\
d<m}}d^{1/2} < \frac{|A|}{2}.
\end{eqnarray*} Hence, $|A\cap \mathbb{Z}_m^*|>|A|/2.$ Denoting
$A\cap \mathbb{Z}_m^*$ again by $A$ we deduce that it suffices to
deal with the case $A\subset \mathbb{Z}_m^*.$

We can also assume that $|A|^3 > 2|AA||A|^2|A+A|/m,$ since otherwise
we are done. Following the proof of
Theorem~\ref{thm:sumprodestPrime}, we obtain
\begin{eqnarray*}
|A|^3\le \frac{2}{m}\sum_{n=1}^{m-1}\Bigl|\sum_{x\in AA}\sum_{a_1\in
A}e_m(nxa_1^{-1})\Bigr|\Bigl|\sum_{a_2\in A}\sum_{y\in
A+A}e_m(n(a_2-y))\Bigr|.
\end{eqnarray*}
For a given divisor $d|m$ we collect together those $n$ for which
$(n,m)=d.$ Thus, denoting $n/d$ by $n$ again, we get
\begin{eqnarray*}
|A|^3 \le \frac{2}{m}\sum_{\substack{d|m\\
d<m}}\sum_{\substack{n=1\\ (n,m/d)=1}}^{m/d}\Bigl|\sum_{x\in
AA}\sum_{a_1\in A}e_{m/d}(nxa_1^{-1})\Bigr|\Bigl|\sum_{a_2\in
A}\sum_{y\in A+A}e_{m/d}(n(a_2-y))\Bigr|.
\end{eqnarray*}
Since $(n,m/d)=1,$ we have
$$
\Bigl(\Bigl|\sum_{x\in AA}\sum_{y\in
A}e_{m/d}(nxa_1^{-1})\Bigr|\Bigr)^2\le
|AA|\sum_{x=0}^{m-1}\Bigl|\sum_{a_1\in
A}e_{m/d}(xa_1^{-1})\Bigr|^2\le dm|AA||A|.
$$
Thus,
$$
|A|^3 \le \frac{2|A|^{1/2}|AA|^{1/2}}{m^{1/2}}\sum_{\substack{d|m\\
d<m}}d^{1/2}\sum_{n=1}^{m/d}\Bigl|\sum_{a_2\in
A}e_{m/d}(na_2)\Bigr|\Bigl|\sum_{y\in A+A}e_{m/d}(ny)\Bigr|.
$$
Using the inequalities
$$
\sum_{n=1}^{m/d}\Bigl|\sum_{a_2\in A}e_{m/d}(na_2)\Bigr|^2\le
m|A|,\qquad \sum_{n=1}^{m/d}\Bigl|\sum_{y\in
A+A}e_{m/d}(ny)\Bigr|^2\le m|A+A|,
$$
we deduce that
$$
|A|^3 \le 2|A|^{1/2}|AA|^{1/2}m^{1/2}|A|^{1/2}|A+A|^{1/2}\sum_{\substack{d|m\\
d<m}}d^{1/2}.
$$
This proves Theorem~\ref{thm:sumprodestCom}.

\bigskip

{\bf Acknowledgement.} The author is grateful to S. V. Konyagin for
very useful remarks. This work was supported by the Project PAPIIT
IN 100307 from UNAM.

\end{document}